\documentclass[12pt]{amsart}
 \usepackage{epsfig}
\usepackage{amssymb}
\usepackage{amsfonts}
\usepackage{amsmath}
\usepackage{array}
\usepackage{rotating}

\newtheorem{example}{Example}[section]

\newtheorem{theorem}[example]{Theorem}

\newtheorem{corollary}[example]{Corollary}

\newtheorem{proposition}[example]{Proposition}

\newtheorem{lemma}[example]{Lemma}

\def\Proof{\noindent \it Proof -- \rm}
\def\qed{\hspace{3.5mm} \hfill \vbox{\hrule height 3pt depth 2 pt width 2mm}
\bigskip}

\def\<{\langle}
\def\>{\rangle}

\def\C{\operatorname{\mathbb C}}
\def\Z{\operatorname{\mathbb Z}}

\def\SG{{\mathfrak S}}

\def\t{{\bf t}}

\def\Peak{{\bf \Pi}}
\def\PP{{\mathcal P}}

\def\Sym{{\bf Sym}}
\def\tensor{\otimes}
\def\Q{\operatorname {Q}}
\def\Des{\operatorname{Des}}
\def\Peak{\operatorname{\sf Peak}}

\def\ch{\operatorname {ch}}
\def\Ind{\operatorname{Ind}}
\def\Res{\operatorname{Res}}

\def\End{\operatorname{End}}
\def\Soc{\operatorname{Soc}}

\def\Tabvrule{\vrule width-0.4pt}       % Difference de largeur
\def\Tabhrule{\hrule \hrule height-0.4pt} % Difference de hauteur
\def\Tabstrut{\vrule height2.2ex % Sur la ligne
                     depth0.8ex  % Sous la ligne
                     width0ex    % centrage horizontal
\relax}

\def\PasCase#1{\omit%
            $\vcenter{\hbox {\vbox to 0.4pt{}}
               \hbox{\makebox[3ex]{\Tabstrut$#1$}}}%
               \Tabvrule$}
\def\PasCasePoint{\PasCase{\cdot}}
\def\DessinCarre#1{%
    \vcenter{\hbox{}\hrule
             \hbox{\vrule\makebox[3ex]{\Tabstrut$#1$}\vrule}\Tabhrule}%
             \Tabvrule}
\def\GenRuban#1{\vcenter{\halign{&$\DessinCarre{##}$\cr#1}}\egroup}

\def\sTabvrule{\vrule width-0.4pt}
\def\sTabhrule{\hrule \hrule height-0.4pt}
\def\sTabstrut{\vrule height1.6ex depth0.6ex width0ex \relax}
\def\sDessinCarre#1{%
    \vcenter{\hbox{}\hrule
             \hbox{\vrule\makebox[2.3ex]%
                  {\sTabstrut$\scriptstyle#1$}\vrule}\sTabhrule}%
             \sTabvrule}
\def\sGenRuban#1{\vcenter{\halign{&$\sDessinCarre{##}$\cr#1}}\egroup}

\def\ruban{%
  \bgroup
  \let\ =\omit
  \let\\=\cr
  \let\x=\times
  \let\.=\PasCasePoint
  \offinterlineskip
  \GenRuban}

\def\sruban{%
  \bgroup
  \let\ =\omit
  \let\x=\times
  \let\\=\cr
  \offinterlineskip
  \sGenRuban}

\title[The peak algebra and the Hecke-Clifford algebras]{The peak algebra and the Hecke-Clifford algebras at $q=0$}

\author[N.~Bergeron, F~Hivert, and J.-Y~Thibon]{Nantel  Bergeron, Florent Hivert, and Jean-Yves Thibon}

\address[Nantel Bergeron]
 {Department of Mathematics and Statistics\\ York University\\
 Toronto, Ontario M3J 1P3\\
 CANADA}
 \email[Nantel Bergeron]{bergeron@mathstat.yorku.ca}
 \urladdr[Nantel Bergeron]{http://www.math.yorku.ca/bergeron}

\address[Florent Hivert and Jean-Yves Thibon]
 {Institut Gaspard Monge, Universit\'e de Marne-la-Vall\'ee \\
5 Boulevard Descartes \\
Champs-sur-Marne \\
77454 Marne-la-Vall\'ee cedex 2\\ 
FRANCE 
}
 \email[Florent Hivert ]{florent.hivert@univ-mlv.fr}
 \email[Jean-Yves Thibon]{jyt@univ-mlv.fr}
 
\date{}

 \thanks{Bergeron is supported in part by CRC, NSERC and PREA}
\begin{document}

\maketitle

\begin{abstract}
Using the formalism of noncommutative symmetric functions, we derive the basic theory of
the peak algebra of symmetric groups and of its graded Hopf dual.  Our main result is to provide
a representation theoretical interpretation of the peak algebra and its graded dual as Grothendieck rings
of the tower of Hecke-Clifford algebras at $q=0$.
\end{abstract}

\section{Introduction}

Studies on the combinatorics of descents in permutations led to the discovery of a pair, $(QSym,\Sym)$, of mutually dual graded Hopf algebras \cite{NCSF1,MR}. Here, $QSym$ is the graded Hopf algebra of quasi-symmetric functions, and its graded dual, $\Sym$, is the graded Hopf algebra of noncommutative symmetric functions. Recent investigations on
the combinatorics of peaks in permutations resulted in the
discovery of an interesting new pair, $(\Peak,\Peak^*)$, of graded Hopf algebras.
The first one, $\Peak$, originally due to Stembridge \cite{St}, 
is a subalgebra of $QSym$. As described in \cite{BMSW},
its graded dual, $\Peak^*$, can therefore be identified as a homomorphic image of $\Sym$. We shall see in the following
that the existence of $\Peak^*$ as well as many of its basic properties
were already implicit in \cite{NCSF2}.

It is known that $\Peak$ can also be obtained as a quotient of $QSym$, in which case
$\Peak^*$ is realized  as a subalgebra of $\Sym$.  On the other hand,
each homogeneous component $\Sym_n$ of $\Sym$ is endowed with another
multiplication, the internal product $*$, such that the resulting algebra is
anti-isomorphic to Solomon's descent algebra of the symmetric group $\SG_n$.
At this stage, a natural question
arises. Is $\Peak_n^*$ stable under this operation? As shown in \cite{Ny},
the answer is yes (it is even a left  ideal of $\Sym_n$), and the corresponding
right ideal of the descent algebra is spanned by the sums of permutations having
a given peak set.
Recent developments  \cite{ABN,ABS,BHW,Scho} unveil many interesting properties and generalizations of  $\Peak$ and $\Peak^*$. Most notably, we find in \cite{ABS} that $\Peak$ is the terminal object in the category of combinatorial Hopf algebras satisfying generalized Dehn-Somerville relations. This reveals some of the significance of $\Peak$ and $\Peak^*$. Our main result demonstrates yet another facet of the importance of these graded Hopf algebras.

We shall start our presentation by showing  that many of  the basic results in the literature related to  $\Peak$ and $\Peak^*$ can be recovered in a very elegant and straightforward way by relying upon the techniques developed
in \cite{NCSF2}. This will be covered in Sections 2,3 and 4.

It is known that the dual pair of Hopf algebras $(QSym,\Sym)$
describes the representation theory of the $0$-Hecke algebras
of type $A$ \cite{NCSF4}. More precisely, $QSym$ and $\Sym$ are
respectively isomorphic to
the direct sums of the Grothendieck groups $G_0(H_n(0))$ and $K_0(H_n(0))$.
We provide here a similar interpretation for the
pair $(\Peak,\Peak^*)$. This is
done by replacing the Hecke algebras with the so-called
Hecke-Clifford algebras, discovered by G. Olshanski \cite{Ol}.
This  new result is presented in Section 5.

We will assume that the reader is familiar with the notation of  \cite{NCSF1,NCSF2}.

 \medskip
 {\footnotesize
 {\it Acknowledgements.-} The computations where done using the package
\texttt{MuPAD-Combinat} by the second author and N.Thi\'ery: \texttt{http://mupad-combinat.sourceforge.net/}. The code can be found in the subdirectory
\texttt{MuPAD-Combinat/lib/experimental/HeckeClifford}. 
 }

\section{The $(1-q)$-transform at $q=-1$}

The main motivation for Stembridge's theory of enriched $P$-partitions,
which led him to the quasi-symmetric peak algebra
\cite{St},  was the study of the quasi-symmetric expansions of Schur's
$Q$-functions \cite{Schur, Mcd}. 
As is well known, these symmetric functions correspond to the
Hall-Littlewood functions with parameter $q=-1$.
The peak algebra is therefore directly related to what we will call 
the  ``$(1-q)$-transform'' at $q=-1$.

In the classical case, the $(1-q)$-transform $\theta_q$ is the ring endomorphism of
$Sym$ defined on the power sums by $\theta_q(p_n)=(1-q^n)p_n$.
In the $\lambda$-ring notation, which is particularly convenient to
deal with such transformations, it reads $f(X)\mapsto {f((1-q)X)}$.
One has to pay attention to the abuse of notation in using the
same minus sign for the $\lambda$-ring and for scalars, though these
operations are quite different. That is, $\theta_{-1}$ maps $p_n$
to $2p_n$ if $n$ is odd, and to $0$ otherwise. Thus, $\theta_{-1}(f(X))=f((1-q)X)_{q=-1}$
is not the same as $f((1+1)X)=f(2X)$.

The main results of \cite{NCSF2} are concerned with the extension
of the $(1-q)$-transform to noncommutative symmetric functions. 
More precisely, consider the abelianization map  $\chi\colon\Sym\to Sym$ which sends the noncommutative alphabet $A$ to the commutative alphabet $X$. We are interested in defining a $(1-q)$-transform on $\Sym$ which commutes with $\chi$.
A consistent
definition of $\theta_q(F)={F((1-q)A)}$ is proposed, and its fundamental
properties are obtained. We briefly recall here the necessary steps. 
One first defines the complete symmetric functions $S_n((1-q)A)$ via their generating series
%%%%%%%%%%%%%%% CHECK
 \cite[Def.~5.1]{NCSF2}
\begin{equation}\sigma_t((1-q)A):=
\sum_{n\ge 0}t^n S_n((1-q)A) = \left(
\sum_{n\ge 0} (-qt)^n\Lambda_n(A) \right)\left(\sum_{n\ge 0}t^n S_n(A)\right)\,,
\end{equation}
and then  $\theta_q$ is defined as  the ring homomorphism such that
$\theta_q(S_n)=S_n((1-q)A)$. Specializing 
%%%%%%%%%%%%%%% CHECK
\cite[Thm.~4.17]{NCSF2} to our case, we obtain 
\begin{equation}\label{internal}
F((1-q)A)=F(A)*\sigma_1((1-q)A),
\end{equation}
where $*$ is the internal product of $\Sym$.

The most important property of $\theta_q$  is
its diagonalization
%%%%%%%%%%%%%%% CHECK
\cite[Thm.~5.14]{NCSF2}:
 there is   a unique family
of Lie idempotents $\pi_n(q)$ (i.e., elements
in the primitive Lie algebra such that $\chi(\pi_n(q))={1\over n}p_n$)
with the property
\begin{equation}\label{proppi}
\theta_q(\pi_n(q))=(1-q^n)\pi_n(q) \,.
\end{equation}
Moreover, $\theta_q$ is semi-simple, and its eigenvalues in $\Sym_n$
are $p_\lambda(1-q)={\prod_i(1-q^{\lambda_i})}$ where $\lambda$
runs over the partitions of $n$. The projectors on the corresponding
eigenspaces are 
 the maps $F\mapsto F*\pi^I(q)$, where for
a composition  $I=(i_1,i_2,\ldots,i_r)$, we let
$\pi^I(q)= \pi_{i_1}(q)\pi_{i_2}(q)\cdots\pi_{i_r}(q)$ 
%%%%%%%%%%%%%%% CHECK
\cite[Sec.~3.4]{NCSF2}.

Another result 
%%%%%%%%%%%%%%% CHECK
\cite[Sec.~5.6.4]{NCSF2},
 which is just a translation of
an important formula due to Blessenohl and Laue \cite{BL}, gives
$\theta_q(R_I)$ in closed form for any ribbon $R_I$. To be more in line with the current literature, we digress slightly from the notation of \cite{NCSF1}. Let $[i,j]=\{i,i+1,i+2,\ldots,j\}$. Given a composition $I=(i_1,i_2,\ldots,i_r)$ of $n$, let $\Des(I)=\{i_1,i_1+i_1,\ldots,i_1+\cdots + i_{r-1}\}\subseteq[1,n-1]$ denote the {\sl descent} set  of $I$. We let $A\Delta B=(A-B)\cup(B-A)$ be the symmetric difference of two sets. Given $A=\{a_1,a_2,\ldots,a_{r-1}\}\subseteq [1,n-1]$ we let $A+1=\{a_1+1,a_2+1,\ldots,a_{r-1}+1\}\subseteq [2,n]$. For a composition $J$ of $n$, one defines $HP(J)=\{a\in \Des(J) \,|\, a\!\ne\! 1, a-1\not\in \Des(J)\}\subseteq [2,n-1]$, and $hl(J)=|HP(J)|+1$. One usually refers to $HP(J)$ as the {\sl peak} set of $J$. We are now in a position to give the formula for $\theta_q(R_I)$ \cite[Lem.~5.38 and Prop.~5.41]{NCSF2}:
\begin{equation}\label{transrub}
R_I((1-q)A)=\sum_{HP(J)\subseteq \Des(I)\Delta(\Des(I)+1)}
(1-q)^{hl(J)}(-q)^{b(I,J)}R_J(A)\,,
\end{equation}
where  $b(I,J)$ is some explicit integer, but is not of any use when $q=-1$.

Setting $q=-1$ in the formulas above leads us
immediately to the peak classes in $\Sym$. We say that a set $P\subseteq[2,n-1]$ is a {\sl peak set} when $a\in P \implies a-1\not\in P$.
For a peak set $P$ let
\begin{equation}\label{defofPi}
\Pi_P=\sum_{HP(I)=P} R_I.
\end{equation}
At $q=-1$, Eq.~(\ref{transrub}) now reads as
\begin{equation} 
\theta_{-1}(R_I)=\sum_{P\subseteq \Des(I)\Delta(\Des(I)+1)}2^{|P|+1}\Pi_P\,,
\end{equation}
which is \cite[Prop.~5.5]{Scho} or \cite[Prop.~5.8]{ABN}.

Let us denote for short $\theta_{-1}$
by a tilde, $\tilde F:= \theta_{-1}(F)$, and let $\widetilde{\Sym}$ be its
image. Since by definition
\begin{equation}
\widetilde{\Sym}=\{F((1-q)A)_{q=-1}\}
\end{equation}
it is immediate that $\widetilde{\Sym}$ is a graded Hopf subalgebra of $\Sym$.
Indeed, $F\mapsto {F((1-q)A)}$ is an algebra morphism, and also a coalgebra
morphism, since \cite[Sec.~5.1]{NCSF2}
\begin{equation}  
\Delta S_n((1-q)A)=\sum_{i+j=n}S_i((1-q)A)\otimes S_j((1-q)A)
\end{equation}
for all values of $q$. Also, it is a left ideal for the internal product,
since by Eq.~(\ref{internal})
\begin{equation} 
\widetilde{\Sym}=\Sym(A)*\sigma_1((1-q)A)_{q=-1}\,.
\end{equation}
We already know that  $\widetilde{\Sym}$ is contained in the subspace $\PP$ of
$\Sym$ spanned by the peak classes. The dimension of the homogeneous component
$\PP_n$ of $\PP$ is easily seen to be equal to the Fibonacci number $f_n$
(with the convention $f_0=f_1=f_2=1$, $f_{n+2}=f_{n+1}+f_n$). 
Indeed, the set of compositions of $n$ having a given peak set $P$ has a unique minimal and maximal element.
The minimal element is a composition where each part except the last is at least $2$.
The maximal element is composed of only 1s and 2s and ends in 1. 
Both sets are obviously of cardinality $f_n$, which is also the number of
compositions of $n$ into odd parts. 

But, thanks to Eq.~(\ref{proppi}), we know that the elements
\begin{equation}\label{oddpibasis}
\pi^I(-1)=\pi_{i_1}(-1)\pi_{i_2}(-1)\cdots \pi_{i_r}(-1)\,,
\end{equation}
where $i=(i_1,\ldots,i_r)$ runs over compositions of $n$ into
odd parts, form a basis of $\widetilde{\Sym}_n$. Hence,
\begin{equation}
\widetilde{\Sym}=\PP_n\,.
\end{equation}
Also, since the commutative image of $\pi_n(q)$ is ${1\over n}p_n$ for
all $q$, this makes it clear that the commutative image
of $\widetilde{\Sym}$ is the subalgebra of $Sym$ generated by
odd power-sums $p_{2k+1}$.

To summarize, we have shown that the peak classes $\Pi_P$ in $\Sym$
form a linear basis of a graded Hopf subalgebra $\PP$ of $\Sym$, which is also
a left ideal for the internal product, and we have described
a basis of it, which is mapped onto products of odd power sums
by the commutative image homomorphism. Since the $\pi_n(-1)$ are
Lie idempotents, this also determines the primitive Lie algebra
of $\PP$ as the free Lie algebra generated by the $\pi_{2k+1}(-1)$.
It is interesting to remark that all this has been obtained without much effort by
setting $q=-1$ in a few formulas of \cite{NCSF2}.

\section{The quasi-symmetric side}

To recover Stembridge's algebra, we have to look at the dual of $\PP$.
Since $\PP$ can be regarded either as a homomorphic image of $\Sym$
(under $\theta_{-1}$) or as a subalgebra of $\Sym$ (spanned by the
peak classes $\Pi_P$), the dual $\PP^*$ can be realized either as a
subalgebra, or as a quotient of $QSym$.

Recall that we have a nondegenerate duality between $QSym$ and $\Sym$
defined by \cite{NCSF1,MR}
\begin{equation}
\<F_I,R_J\>=\delta_{IJ}\,.
\end{equation}

Let us first consider the noncommutative peak algebra $\PP=\Peak^*$ as the image
of the Hopf epimorphism $\varphi=\theta_{-1}$. Then, the adjoint map
\begin{equation}
\varphi^*:\ \PP^* \longrightarrow QSym
\end{equation}
is an embedding of Hopf algebras. The duality between $\PP$ and $\PP^*$
is given by
\begin{equation}
\<\varphi(F),G\>=\<F,\varphi^*(G)\>\,.
\end{equation}
Hence, if we denote by $\Pi_P^*$ the dual basis of $\Pi_P$, we have for
any ribbon $R_I$ with descent set $D=\Des(I)$
\begin{equation}
\<\varphi^*(\Pi_P^*),R_I\>=\<\Pi_P^*,\varphi(R_I)\>
=\begin{cases} 2^{|P|+1} &\text{if $P\subseteq D\Delta(D+1)$}\\ 0 &\text{otherwise}\,.
 \end{cases}
\end{equation}
Thus, in its realization as a subalgebra of $QSym$, $\PP^*$ is spanned by
Stembridge's quasi-symmetric functions
\begin{equation}
\Theta_P= \varphi^*(\Pi_P^*)=
2^{|P|+1}\sum_{P\subseteq\Des(I)\Delta(\Des(I)+1)}F_I\,.
\end{equation}

Note also that thanks to the identity $(1+q)(1-q)=1-q^2$, the kernel of
$\varphi$ is seen to be the ideal of $\Sym$ generated by the $S_n((1+q)A)_{q=-1}$
for $n\ge 1$. These are the $\chi_n$ of \cite{BHW}.

Finally, we can also consider $\PP$ as an abstract algebra with basis
$(\Pi_P)$, and define a monomorphism $\psi:\ \PP\rightarrow \Sym$
by
\begin{equation}
\psi(\Pi_P)=\sum_{HP(I)=P}R_I\,.
\end{equation}
Then, its adjoint $\psi^*:\ QSym\rightarrow \PP^*$ is an epimorphism.
The product map $\vartheta=\varphi^*\circ\psi^*: \PP^*\rightarrow \PP^*$
has been considered by Stembridge \cite{St}, and its diagonalization is
given in \cite{BHW}. We can easily recover its properties from the results
of the previous section, since clearly $\vartheta=(\psi\circ\varphi)^*$ coincides
with $\theta_{-1}$. Its eigenvalues are then 
the integers $2^{\ell(\lambda)}$,
where $\lambda$ runs over partitions into odd parts. 
The spectral projectors are again constructed from the idempotents
$\pi_\lambda(-1)$. Precisely, the projector onto the eigenspace
associated with the eigenvalue $2^k$ of $ \vartheta$ in $QSym_n$ is the adjoint of
the endomorphism of $\Sym_n$ given by $F\mapsto F*U_k$ where
$U_k=\sum\pi_\lambda(-1)$, the sum being over all odd partitions of $n$
with exactly $k$ parts. The dimensions of these eigenspaces can also be easily
computed.

\section{Miscellaneous related results}

Here are some more results related to the recent literature. We choose to include them here for completeness. 

\subsection{Noncommutative tangent numbers}

By definition, $\widetilde{\Sym}$ is generated by the $\tilde S_n$, $n\ge 1$.
If we set $q=-1$ in \cite[Prop. 5.2]{NCSF2} we establish that  $\tilde S_n=2H_n$ for $n\ge 1$, where $H_0=1$ and
\begin{equation}
H_n=\sum_{k=0}^{n-1}R_{1^k,n-k}\,.
\end{equation}
Then \cite[Prop. 5.24]{NCSF1} gives us that
\begin{equation}\label{th} 
H=\sum_{n\ge 0}H_n=(1-\t)^{-1}
\end{equation}
where $\t$ is the (left) noncommutative hyperbolic tangent
\begin{equation}   
\t={\rm TH} = \sum_{k\ge 0}(-1)^kT_{2k+1}\,,
\quad
T_{2k+1}=R_{2^k1}\,.
\end{equation}
Hence,  $\widetilde{\Sym}$ is contained in the subalgebra generated by
the $T_{2k+1}$, and since we already know that the dimension of
 $\widetilde{\Sym}_n$ is the number of odd compositions of $n$, we have
in fact equality. Thus, the $T^I=T_{i_1}\cdots T_{i_r}$ ($I$ odd)
form a multiplicative basis of  $\widetilde{\Sym}$ (this is the same
as the basis $\Gamma^P$ of \cite{Scho}).

\subsection{Peak Lie idempotents}

We have
seen in Section 2, Eq.~(\ref{oddpibasis}), that the $\pi^I(-1)$, $I$ odd, form a basis of  $\widetilde{\Sym}$,
so that  $\widetilde{\Sym}_n$ contains Lie idempotents iff $n$ is
odd. In \cite{Scho}, the images $\tilde L_n$ of some classical Lie idempotents 
$L_n$ are calculated. 

The images $\tilde\Psi_n=\Psi_n((1-q)A)$ of the Dynkin elements $\Psi_n$ are given in closed
form for any $q$ in \cite[Prop. 5.34]{NCSF2}. It suffices to set $q=-1$
in this formula to obtain \cite[Prop. 7.3]{Scho}.

On the other hand \cite[Prop. 7.2]{Scho} gives an interesting new
formula for $\tilde\Phi_n$. Yet, the first part of the analysis can be simplified by applying
Eq.~(\ref{th}) to the calculation of the generating
series $\log \tilde\sigma_1$. Indeed,
\begin{eqnarray*}
\log \tilde\sigma_1 &=& \log(1+\t)-\log(1-\t) \\
&=& 2\sum_{k\ge 0}{\t^{2k+1}\over 2k+1}\\
&=&  2\sum_{I \ {\rm and}\ \ell(I)\ {\rm odd}}
{(-1)^{(|I|-\ell(I))/2} \over \ell(I)}T^I\,.
\end{eqnarray*}

Finally, to obtain the image of Klyachko's idempotent $\tilde K_n(q)$ (see  \cite[Prop.~6.3]{NCSF2} for a definition of $K_n(q)$) one has to set $t=-1$ in
 \cite[Prop.~8.2]{NCSF2}.

\subsection{Structure of the Peak algebras $(\PP_n,*)$}

Using the construction of \cite[Sec.~3.4]{NCSF2} restricted to odd partitions $\lambda$ of $n$, it follows from Eq.~(\ref{proppi}) that the idempotents $E_\lambda(\pi(-1))$,
associated to the sequence $\pi_n(-1)$, form a complete
set of orthogonal idempotents of $\PP_n$. Regarding $\PP_n$ as a quotient
of the descent algebra makes it clear that the left ideals $\PP_n*E_\lambda(\pi(-1))$
are the indecomposable projectives modules of $\PP_n$. We obtain explicitly the multiplicative
structure of $(\PP_n,*)$ by adapting 
\cite[Lem.~3.10]{NCSF2} to the sequence $\pi_n(-1)$ (instead of $\Psi_n$), and then imitating
the rest of the argument presented there for the descent algebra.

\subsection{Hall-Littlewood basis}

The peak algebra $\PP=\Peak^*$ can be regarded as a noncommutative version
of the subalgebra of $Sym$ spanned by the Hall-Littlewood functions
$Q_\lambda(X;-1)$, where $\lambda$ runs over strict partitions.
Actually, it is easy to show that the noncommutative Hall-Littlewood
functions of \cite{BZ,Hi} at $q=-1$ yield two different analogous bases of
$\PP$. We do it here for \cite{Hi} but a similar argument can be applied to \cite{BZ}.

Recall that the polynomials $H_I(A;q)$ of \cite{Hi} are defined as
noncommutative analogues of the $Q'_\mu=Q_\mu(X/(1-q);q)$. To obtain
the correct analogues of Schur's $q$-functions, one has to apply
the $(1-q)$ transform before setting $q=-1$.

\begin{proposition}
The specialized noncommutative Hall-Littlewood functions
\begin{equation}
\Q_I=H_J((1-q)A;q)_{q=-1},
\end{equation}
where $I$ runs over all peak compositions, form a basis of $\PP$.
\end{proposition}

Indeed, the factorization of $H$-functions at roots of unity imply that
\begin{equation}
\Q_I=\Q_{i_1i_2}\Q_{i_3i_4}\cdots Q_{i_{2k-1}i_{2k}}\Q_{i_{2k+1}}
\end{equation}
(where $i_{2k+1}=0$ if $I$ is of even length), and simple calculations
yield
\begin{itemize}
\item $\Q_n=2\Pi_\emptyset$,
\item $\Q_{n-1,1}=2(\Pi_{\{n-1\}}+\Pi_\emptyset)$,
\item and for $2\le k\le n-2$,
    $\Q_{k,n-k}=4(\Pi_{\{k\}}+\Pi_{\{k+1\}}+\Pi_\emptyset)$,
\end{itemize}
where $\Pi_P$ is defined in Eq.~(\ref{defofPi}).
From this, it is straightforward to prove that the family $\Q_I$
is triangular with respect to the family $\Pi_P$, and hence
the proposition follows.

%%%%%%%%%%%%%%%%%%%%%%%%%%%%%%%%%

\section{Representation theory of the $0$-Hecke-Clifford algebras}

The character theory of symmetric groups (in characteristic 0),
as worked out by Frobenius, can be summarized as follows.
Let  $R_n$ denote the free abelian group spanned by
isomorphism classes of irreducible representations of $\C\SG_n$.
Endow the direct sum
\begin{equation}
R=\bigoplus_{n\ge 0} R_n
\end{equation}
with the addition corresponding to direct sum, and 
multiplication $R_m\otimes R_n\rightarrow R_{m+n}$ corresponding
to induction from $\SG_m\times\SG_n$ to $\SG_{m+n}$ via the natural
embedding. The linear map sending the class of an irreducible
representation $[\lambda]$ to the Schur function $s_\lambda$
is then a ring isomorphism between $R$ and $Sym$ (see, e.g., \cite{Mcd}). Moreover, we can define a structure of  graded Hopf algebra on $R$ with comultiplication corresponding to restrictions from $\SG_{n}$ down to $\SG_k\times\SG_{n-k}$ and summing over $k$. The linear map above gives rise to an isomorphism of graded Hopf algebras.

It is known  that the pair of graded Hopf algebras
$(\Sym, QSym)$ admits a similar interpretation, in terms
of the tower of
the $0$-Hecke algebras $H_n(0)$ of type $A_{n-1}$ (see \cite{NCSF4}). 
Recall that the (Iwahori-) Hecke algebra $H_n(q)$ is the $\C$-algebra generated by
elements
$T_i$ for $i<n$ with the relations:
\begin{alignat}{2}
  T_i^2&=(q-1)T_i+q    \qquad&&\text{for $1\leq i \leq n-1$,}\notag\\
  T_iT_j&=T_jT_i \qquad&&\text{for $|i-j|>1$,}\\
  T_iT_{i+1}T_i &=T_{i+1}T_iT_{i+1} \qquad&&\text{for $1\leq i \leq n-2$,}\notag
\end{alignat}
(here, we assume that $q\in \C$).
The $0$-Hecke algebra is obtained by setting $q=0$ in these
relations. Then, the first
relation becomes $T_i^2=-T_i$. If we denote by $G_n=G_0(H_n(0))$
the Grothendieck group of the category of finite dimensional $H_n(0)$-modules, and by $K_n=K_0(H_n(0))$ the Grothendieck group of the
category of projective $H_n(0)$-modules,
the direct sums ${\mathcal G}=\bigoplus_{n\ge 0}G_n$ and
${\mathcal K}=\bigoplus_{n\ge 0}K_n$, endowed with the same operations
as above, are respectively isomorphic with $QSym$ and $\Sym$. 

The aim of this final section is our main result: to provide a similar interpretation
for the pair $(\PP,\PP^*)$. The relevant tower of algebras
is the $0$-Hecke-Clifford algebras, which are degenerate
versions of Olshanski's Hecke-Clifford algebras.

\subsection{Hecke-Clifford algebra}
%%%%%%%%%%%%%%%%%%%%%%%%%%%%%%%%%%%%%%%%%%%%%%%%%%%%%%%%%%

The complex Clifford algebra $Cl_n$ is generated by $n$ elements 
$c_i$ for $i\leq n$ with the
relations
\begin{equation}
  c_ic_j = - c_jc_i \quad\text{for $i\neq j$\qquad and \qquad}
  c_i^2 = -1.
\end{equation}
For each subset $D = \{i_1 < i_2 < \dots i_k \}\subset\{1 \dots n\}$, we
denote by $c_D$ the product
\begin{equation}
  c_D := \prod_{i\in D}^\rightarrow c_i = c_{i_1} c_{i_2} \dots c_{i_k}\,.
\end{equation}
It is easy to see that $(c_D)_{D\subset\{1 \dots n\}}$ is a basis of the
Clifford algebra.

The Hecke-Clifford superalgebra \cite{Ol} is the unital $\C$-algebra
generated by the $c_i$, and $n-1$ elements $t_i$ satifying
the Hecke relations in the form
\begin{alignat}{2}
t_i^2&=(q-q^{-1})t_i+1    \qquad&&\text{for $1\leq i \leq n-1$,}\notag\\
t_it_j&=t_jt_i \qquad&&\text{for $|i-j|>1$,}\\
t_it_{i+1}t_i &=t_{i+1}t_it_{i+1} \qquad&&\text{for $1\leq i \leq n-2$.}\notag
\end{alignat} 
and the cross-relations
\begin{alignat}{2}
t_ic_j&= c_jt_i                \qquad&&\text{for $i\neq j,j+1$,}\notag\\
t_ic_i&= c_{i+1}t_i            \qquad&&\text{for $1\leq i \leq n-1$,}\\
(t_i+q^{-1})c_{i+1} &= c_{i}(t_i-q) \qquad&&\text{for $1\leq i \leq n-1$.}\notag
\end{alignat}       
The Hecke-Clifford algebra has a natural  $\Z_2$-grading, for which the $t_i$ are
even and the $c_j$ are odd. Henceforth, it will be considered as a superalgebra.
                                                              
Setting $t_i=q^{-1}T_i$ and taking the limit $q\rightarrow 0$ after
clearing the denominators, we obtain the $0$-Hecke-Clifford
algebra $HCl_n(0)$, which is generated by the $0$-Hecke algebra and the
Clifford algebra, with the cross-relations
\begin{alignat}{2}
  T_ic_j&= c_jT_i                \qquad&&\text{for $i\neq j,j+1$,}\notag\\
  T_ic_i&= c_{i+1}T_i            \qquad&&\text{for $1\leq i \leq n-1$,}
\label{crossrel}  \\
  (T_i+1)c_{i+1} &= c_{i}(T_i+1) \qquad&&\text{for $1\leq i \leq n-1$.}\notag
\end{alignat}
Let $\sigma = \sigma_{i_1}\dots\sigma_{i_p}$ be a reduced word
for a permutation $\sigma\in\SG_n$. The defining
relations of $H_n(q)$ ensure that the element $T_{\sigma} :=
T_{i_1}\dots T_{i_p}$ is independent of the chosen reduced word for $\sigma$.
The family $(T_{\sigma})_{\sigma\in\SG_n}$ is a basis of the Hecke
algebra.              
Thus a basis for $HCl_n(q)$ is given by
$(c_D T_\sigma)_{D\subset \{1,\dots, n\}, \sigma \in\SG_n}$,
and consequently, the dimension of $HCl_n(q)$ is $2^n\,n!$ for all $q$. 

\subsection{Quasi-symmetric Characters of induced modules}
%%%%%%%%%%%%%%%%%%%%%%%%%%%%%%%%%%%%%%%%%%%%%%%%%%%%%%%%%%

Since  $H_n(0)$ is  the sub-algebra of $HCl_n(0)$ generated
by the $T_i$, our main tool  in the sequel will
be  the induction process with respect to  this
inclusion. Let us recall some known facts about the representation
theory of $H_n(0)$. There are $2^{n-1}$ simple $H_n(0)$-modules.
These are all one dimensional and can be conveniently labelled by
compositions $I$ of $n$. The  structure of the simple module $S_I :=
\C \epsilon_I$ is given by
\begin{equation}\label{actionTj}
  T_j \epsilon_I = 
  \begin{cases}
    -\epsilon_I &\text{if $j \in \Des(I)$,}\\
    0           &\text{otherwise.}\\
  \end{cases}
\end{equation}
As described in \cite{NCSF4},  there is an isomorphism $\ch\colon  {\mathcal G}\to QSym$ which we call the {\sl Frobenius characteristic}. 
This maps the simple module $S_I$ to the quasi-symmetric function $F_I$. 

Let us define the $HCl_n(0)$-module $M_I$ as the module induced by $S_I$
through the natural inclusion map, that is
\begin{equation}
M_I := \Ind_{H_n(0)}^{HCl_n(0)} (S_I) = HCl_n(0) \tensor_{H_n(0)} S_I\,.
\end{equation}
A basis for $M_I$ is given by $(c_D \epsilon_I)_{D\subset \{1,\dots, n\}}$. A
basis element can be depicted conveniently as follows. The boxes of the ribbon
diagram associated with $I$ are numbered from left to right 
and from top to bottom. 
We
put a ``$\times$" in the $i$-th box if $i\in D$. For example
$c_{\{1,3,4,6\}}\,\epsilon_{(2,1,3)} = c_1c_3c_4c_6\,\epsilon_{(2,1,3)}$ is
depicted by
\begin{equation}\label{ribexam}
  (2,1,3) = \ruban{ 1  & 2  \\
                    \  & 3  \\
                    \  & 4  & 5  & 6\\}
        \hskip2cm
        {c_{\{1,3,4,6\}}\,\epsilon_{(2,1,3)}} = 
  \ruban{ \x &    \\
          \  & \x \\
          \  & \x &    & \x \\}
\end{equation}
We can graphically view the set $\Des(I)$ as the set of boxes with a box below, and the set $HP(I)$ as the set of boxes with  boxes below and  to the left. In the example above, $\Des(2,1,3)=\{2,3\}$ are the boxes labeled $2$ and $3$, and $HP(2,1,3)=\{2\}$ is only the box $2$.

We now remark that  $T_i$ acts only on the $i$-th and $i+1$-st boxes. On the graphical
representation, drawing only the boxes $i$ and $i+1$, the  rules~(\ref{crossrel}) read
\def\rub#1{\text{$\sruban{#1}$}}
\begin{equation}
\label{equ.commRelat}
  \begin{array}{rcrcrcl}
    T_i\ \rub{   &   \\} &=& 0 \ \null& \qquad\qquad\qquad&
    T_i\ \rub{  &  \x\\} &=& -\rub{ & \x\\} + \rub{\x&   \\}
    \\[2mm]
    T_i\ \rub{\x&    \\} &=& 0 \ \null&&
    T_i\ \rub{\x&  \x\\} &=& -\rub{\x& \x\\} + \rub{  &   \\} 
    \\[2mm]
    T_i\ \rub{\x\\   \\} &=& -\rub{  \\ \x\\} &&
    T_i\ \rub{  \\   \\} &=& -\rub{  \\   \\} 
    \\[4mm]
    T_i\ \rub{\x\\ \x\\} &=& \  \rub{  \\   \\} &    &
    T_i\ \rub{  \\ \x\\} &=& -\rub{  \\ \x\\} 
  \end{array}
\end{equation}

At this point, we can make a couple of useful remarks. Looking at the support of the relation (33), we define
%%%%%%%%%%%%%%
\begin{equation}
  \begin{array}{c@{\qquad}c@{\qquad}c@{\qquad}c}
    \rub{  &  \x\\} \rightarrow \rub{\x&   \\}\,,
    &
    \rub{\x&  \x\\} \rightarrow \rub{  &   \\}\,,
    &
    \rub{\x\\   \\} \rightarrow \rub{  \\ \x\\}\,, 
    &
    \rub{\x\\ \x\\} \rightarrow \rub{  \\   \\}\,.
  \end{array}
\end{equation}
  These relations can be interpreted as the cover relation 
 of a (partial) 
order $\leq_I$ on the subset $D$ of
$\{1,\dots, n\}$. Here is a picture of the Hasse diagram of this order for the
composition $I=(2,1,1)$. The poset clearly has two components corresponding to
the two $\Z_2$-graded homogenous components of $M_{(211)}$.

\vskip0.7cm
{\begin{center}
\noindent\epsfig{file=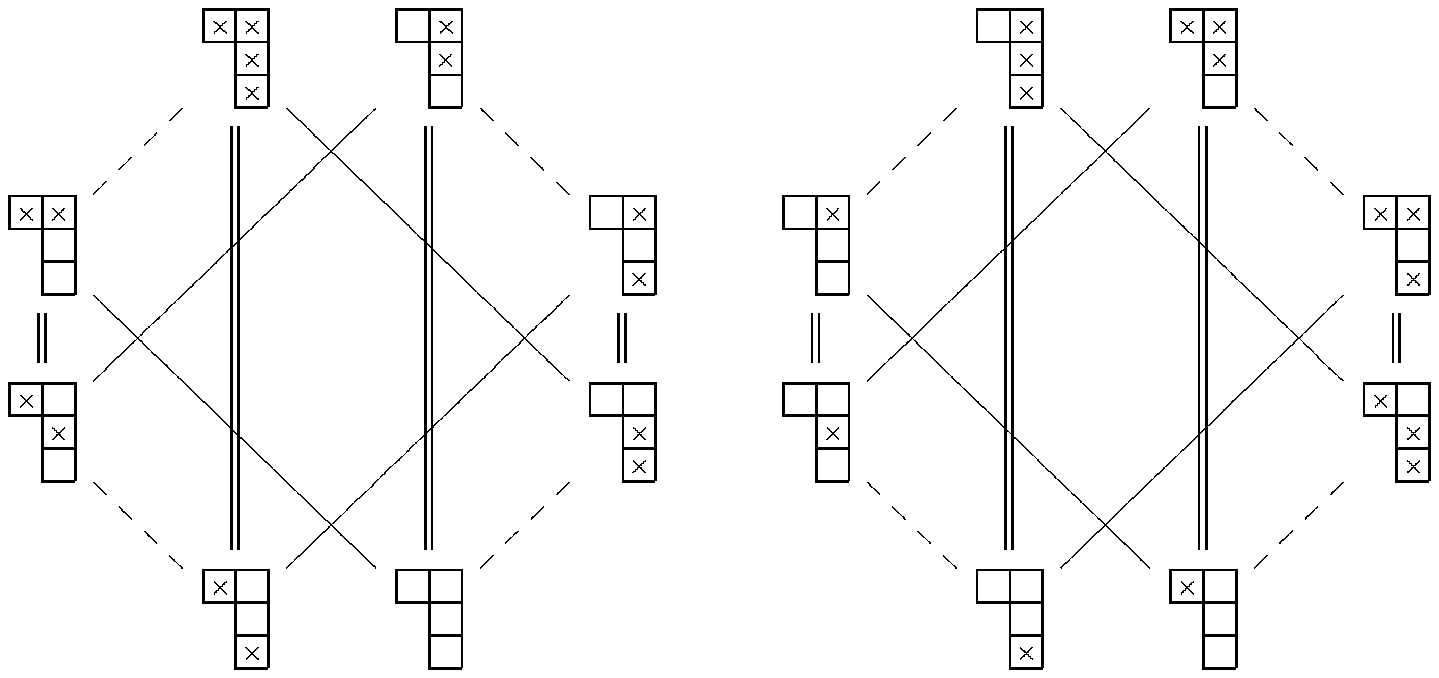,width=\textwidth}

\smallskip
  The Hasse diagram of the order $\leq_{(2,1,1)}$. 
\end{center}
\medskip}
The importance of this order comes from the following lemma, a direct consequence of Eq.~(\ref{equ.commRelat}).
\begin{lemma}
  \label{lem.triangleAction}
  The action of each $T_i$ is triangular with respect to the order $\leq_I$,
  that is for all $D$,
\begin{equation}
  T_i\ c_D \epsilon_I = \alpha(i,I,D)\ c_D \epsilon_I +
  \text{smaller terms}
\end{equation}
with $\alpha(i,I,D)\in\{0,-1\}$.
\end{lemma}

The $\alpha(i,I,D)$ are the eigenvalues of the $T_i$. They are equal to
$0$ in the leftmost columns of Eq.~(\ref{equ.commRelat}) and to $-1$ in
the rightmost ones.

A second consequence of Eq.~(\ref{equ.commRelat}) is that in a vertical (two boxes) diagram,
 the eigenvalue
depends only on the content of the upper box whereas in a horizontal diagram it
depends only on the content of the rightmost box. Thus the content of the boxes 
without a box below  or on the left does not matter for computing the eigenvalues. 
This can be translated into the following lemma.
\begin{lemma}
  \label{lem.defvalley}
  Suppose that $k\in\{1\dots n\}$ is such that  $k$ is not a descent of $I$
  and has no box to its left. 
Then for all $D$, the eigenvalues
  $\alpha(i,I,D)$ satisfy
  \begin{equation}
    \alpha(i,I,D) = \alpha(i,I,\, \{k \} \cup D).
  \end{equation}
  Such a $k$ is called a \emph{valley} of the composition $I$.
\end{lemma}

 Note that  $1$ and $n$ can be valleys. There is
obviously one more valley than the number of peaks.

Thanks to the order $\leq_I$, one can easily describe the structure of the
restriction of $M_I$ to $H_n(0)$. Our first goal is to get a composition
series of $\Res_{H_n(0)} M_I$ in order to compute its Frobenius
characteristic. This can be done as follows. 
Let us choose a linear extension $D_1, D_2 \dots D_{2^n}$
of $\leq_I$. 
For $k\ge 1$, define
\begin{equation}
  M_I^k = \bigoplus_{l\leq k}\ \C \, c_{D_k}\, \epsilon_I,
\end{equation}
and $M_I^0 := \{0\}$. Then, thanks to Lemma \ref{lem.triangleAction},
$M_I^k$ is clearly a sub-module of $\Res_{H_n(0)} M_I$, and
\begin{equation}
  \{0\} = M_I^0 \subset M_I^1 \subset M_I^2 \subset \dots
  \subset M_I^{2^n} = \Res_{H_n(0)} M_I
\end{equation}
is a composition series of the module $\Res_{H_n(0)} M_I$. 
Let us compute the simple
composition factors of the module $S_{D_i,I}=\C\epsilon_{K(D_i,I)} := M_I^i / M_I^{i-1}$. For $1\le k<n$, the generator
$T_k$ acts as $T_k \epsilon_{K(D_i,I)}=\alpha(k, I, D_i)\epsilon_{K(D_i,I)}$. The eigenvalue $\alpha(k, I, D_i)$ equals  $-1$ if 
\begin{equation}
(k+1\in D_i \text{\ and\ } k\not\in \Des(I)) 
\quad\text{or}\quad
(k  \not\in D_i \text{\ and\ } k     \in \Des(I)),
\end{equation}
and $0$ otherwise.
Hence, according to Eq.~(\ref{actionTj}), $\ch(S_{D_i,I})=F_{K}$ where $\Des(K)=\Des(K(D_i,I))=\{1\le k<n\, |\, \alpha(k, I, D_i)=-1\}$.
When $p$ is a peak of $I$, that is $p\ne1$,
 $p-1\not\in\Des(I)$ and $p\in\Des(I)$, then
$\big|\{p-1,p\}\cap \Des(K)\big|=1$. Indeed, if
$p\in D_i$ then $p-1\in \Des(K)$ and $p\not\in \Des(K)$ and if $p\not\in D_i$ then $p-1\not\in \Des(K)$ and $p\in \Des(K)$.
Thus, $P=HP(I)\subseteq\Des(K)\Delta(\Des(K)+1)$. Moreover, for $k\not\in P\cup(P-1)$, we can always find a $D_i$ such that $k\in \Des(K(D_i,I))$.  All $K$ such that $P\subseteq\Des(K)\Delta(\Des(K)+1)$
can be obtained,
 and thanks to Lemma \ref{lem.defvalley} there are
$2^{|P|+1}$ sets $D_i$ giving the same $F_K$. Thus, we have proved the following
proposition. 

\begin{proposition}\label{Fcharac}
  The Frobenius characteristic of $\Res_{H_n(0)} M_I$ depends only on the peak
  set $P$ of the composition $I$ and is given by Stembridge's
 $\Theta$ function 
  \begin{equation}
    \ch(\Res_{H_n(0)} M_I) = 
    \Theta_P = 
    2^{|P|+1}\sum_{P\subseteq\Des(K)\Delta(\Des(K)+1)}F_K\,.
    \end{equation}
\end{proposition}
\qed

\subsection{Homomorphisms between  induced modules}
%%%%%%%%%%%%%%%%%%%%%%%%%%%%%%%%%%%%%%%%%%%%%

The previous proposition suggests that $\Res_{H_n(0)} M_I$ is 
isomorphic to $\Res_{H_n(0)} M_J$ iff
$I$ and $J$ have the same peak sets.
This is actually true, and in fact,  $M_I$ and
$M_J$ are even isomorphic as $HCl_n(0)$-supermodules,
as  we will establish now.

\begin{theorem}\label{MIEndo}
  Let $I$ be a composition with valley set $V$, and  let $Cl_V$ be the
  subalgebra of $Cl_n$ generated by $(c_v)_{v\in V}$. For $c \in Cl_V$ define
  a map $f_c$ from $M_I$ to itself by
  \begin{equation}
   f_c(x \epsilon_I) = x c \epsilon_I
   \qquad
   \text{for all $x\in Cl_n$}.
  \end{equation}
  Then $c \mapsto f_c$ defines a right action of $Cl_V$ on $M_I$ which commutes
  with the left $HCl_n(0)$-action. Moreover the map $c \mapsto f_c$ is a
   graded isomorphism from $Cl_V$ to $\End_{HCl_n(0)}(M_I)$.
\end{theorem}

\Proof 
Since $M_I$ is freely generated as a $Cl_n$-module by $\epsilon_I$, a morphism $f\in\End_{HCl_n(0)}(M_I)$ is determined by $f(\epsilon_I)=x\epsilon_I$ for $x\in Cl_n$. 
On the other hand, for $x\in Cl_n$, a map $f_x(\epsilon_I)=x\epsilon_I$ is in $\End_{HCl_n(0)}(M_I)$ if and only if $T_j f_x(\epsilon_I) = T_j x\epsilon_I = xT_j \epsilon_I$ for all $1\le j<n$. Thus, to prove the theorem, it is
sufficient to see that $f_x\in\End_{HCl_n(0)}(M_I)$ if and only if $x\in Cl_V$.
Equivalently,
\begin{equation}\label{conditionPeak}
  \text{for all $1\le j<n$}
  \qquad\quad
T_j x \epsilon_I =
  \begin{cases}
    -x \epsilon_I &\text{if $j \in \Des(I)$,}\\
    0             &\text{otherwise,}\\
  \end{cases}
\end{equation}
 if and only if $x\in Cl_V$.  

Let us first assume that $x=c_D$ for $D\subseteq V$. This means that in the graphical representation of $x\epsilon_I$ there is no box below nor to the left of a box with a  ``$\times$''. If $j\in \Des(I)$, the lower two equations of the right column of Eq.~(\ref{equ.commRelat}) then
show that $T_jx\epsilon_I=-x\epsilon_I$. The top two equations on the left show that if $j\not\in \Des(I)$ then $T_jx\epsilon_I=0$. By linearity, we get that if $x\in Cl_V$ then Eq.~(\ref{conditionPeak}) holds.  Conversely, let $x\in Cl_n$ satisfy Eq.~(\ref{conditionPeak}).  Let $c_D\epsilon_I$ be in the support of $x$, minimal with respect to $\le_I$.  If $D\not\subseteq V$ then there is a box $j$ with a ``$\times$'' and a box below or to the left.
Using Lemma~\ref{lem.triangleAction} and Eq.~(\ref{equ.commRelat}) this would be a contradiction to Eq.~(\ref{conditionPeak}). Hence $c_D\in Cl_V$. We can subtract it from $x$ and repeat the argument above recursively to conclude that $x\in Cl_V$.
  \qed

\begin{theorem}
  The induced supermodules $M_I$ and $M_J$ are  isomorphic if and only if the
  peak sets of $I$ and $J$ coincide. 
\end{theorem}

\proof
One direction of this theorem is implied by the previous section.
If $M_I$ is isomorphic to $M_J$, then we must have that $\Res_{H_n(0)} M_I$ is isomorphic to $\Res_{H_n(0)} M_J$. In particular, they must have the same Frobenius characteristic. Thanks to Prop.~\ref{Fcharac}, $\ch(\Res_{H_n(0)} M_I)$ depends only on the peak set of $I$. Thus,
if $M_I$ and $M_J$ are isomorphic then $I$ and $J$ have the same peak sets

The converse will follow once we construct explicit isomorphisms between  any modules
$M_I$ and $M_J$ in the same peaks class ($HP(I)=HP(J)$), such that $I$ and $J$ differ exactly by one descent. Isomorphisms between any modules
$M_I$ and $M_J$ in the same peaks class in general will be obtained by composition of the constructed ones.

Let $I=J\cap\{k\}$ be such that $HP(I)=HP(J)$. Graphically, there are two possible cases to consider:
\begin{equation}
\label{samepeak1}
J=\begin{array}{ccc}
     \ddots &&\\
     &\rub{  \\  & \\}&\\
     && \ddots\\
   \end{array} \qquad\qquad\quad
I=\begin{array}{ccc}
     \ddots &&\\
     &\rub{  \\   \\ \\}&\\
     && \ddots\\
   \end{array}
\end{equation}
or
\begin{equation}
\label{samepeak2}
J=\begin{array}{ccc}
     \ddots &&\\
     &\rub{  &  & \\}&\\
     && \ddots\\
   \end{array} \qquad\qquad\quad
I=\begin{array}{ccc}
     \ddots &&\\
     &\rub{  \\   & \\}&\\
     && \ddots\\
   \end{array}
\end{equation}
In the case (\ref{samepeak1}), we construct a map $f$ which sends $\epsilon_I\mapsto \eta=(c_{\{k,k+1\}}-1)\epsilon_J$. We remark that both $\eta$ and $\epsilon_I$ are even.  Furthermore, 
\begin{equation}\label{conditioneta}
  \text{for all $1\le i<n$}
  \qquad\quad
T_i \eta =
  \begin{cases}
    -\eta&\text{if $i \in \Des(I)$,}\\
    0             &\text{otherwise.}\\
  \end{cases}
\end{equation}
Indeed, for $i\not\in\{k-1,k\}$, the $T_i$ commute with $c_{\{k,k+1\}}$. Moreover  $i\in \Des(I)$ if and only if $i\in\Des(J)$, hence the Eq.~(\ref{conditioneta}) follows in these cases. If $i=k-1\in\Des(I)$, then $T_{k-1}\eta=(-c_{\{k-1,k+1\}}T_k + c_{\{k-1,k+1\}}-c_{\{k,k+1\}}-T_k)\epsilon_J=-\eta$, and if $i=k\in\Des(I)$, then $T_{k}\eta=(-c_{\{k,k+1\}}T_k - c_{\{k,k+1\}}+1-T_k)\epsilon_J=-\eta$.
This allows us to define a non-trivial $HCl_n(0)$ supermorphism $f\colon M_I\to M_J$ where
$f(c_D\epsilon_I)=c_D\eta$. Thanks to Eq.~(\ref{conditioneta}), the submodule spanned by $\eta$ in $M_J$ is isomorphic to $M_I$. But since both spaces have the same dimension we have that $f$ is an isomorphism.

For the case (\ref{samepeak2}), we proceed in the same way, sending $\epsilon_J\mapsto (c_{\{k,k+1\}}+1)\epsilon_I$. This constructs a graded isomorphism from $M_J$ to $M_I$.
\qed

\subsection{Simples supermodules of $HCl_n(0)$}
%%%%%%%%%%%%%%%%%%%%%%%%%%%%%%%%%%%%%%%%%%%%  

We are now in a position to construct the simple supermodules of $HCl_n(0)$.
Our approach is similar to Jones and Nazarov \cite{JN}.
Let $I$ be a composition with peak set  $P$ and valley set
$V=\{v_1,v_2,\dots v_k\}$. 
Choose a minimal even idempotent of the Clifford
superalgebra $Cl_V$. For example
\begin{equation}
  \label{eq.minImdemCL}
  e_I :=  {1\over 2^l}
  (1 + \sqrt{-1} c_{v_1}c_{v_2})
  (1 + \sqrt{-1} c_{v_3}c_{v_4})
  \dots
  (1 + \sqrt{-1} c_{v_{2l}}c_{v_{2l+1}})\,,
\end{equation}
where $l := \lfloor \frac{k}2 \rfloor = \lfloor \frac{|P|+1}2 \rfloor$.
Define $HClS_I :=Cl_n e_I \epsilon_I$ as the $HCl_n(0)$-module generated by
$e_I\epsilon_I$. One has to show that $HClS_I$ does not depend on the
chosen minimal idempotent $e_I$, but this is an easy consequence of the representation theory of $Cl_V$
which is know to be supersimple (see e.g. \cite{Jo}).  Suppose that $e_I$ and $e'_I$ are two
minimal even idempotents of $Cl_V$.  Since $Cl_V$, is supersimple, there
exist $x,y,x',y'\in Cl_V$ such that $e_I=x'e'_I y'$ and $e'_I = xe_Iy$. Then
by the representation theory of $Cl_V$, we know that $f_y$ and $f_{y'}$ are two
mutually reciprocal isomorphisms between $Cl_V e_I$ and $Cl_V e'_I$ and hence
between $Cl_n e_I \epsilon_I$ and $Cl_n e'_I \epsilon_I$.

When $n$ is even $Cl_V e_I$ has dimension $2^{\frac{|P|+1}2}$ and when $n$ is odd the dimension is $2^{\frac{|P|+1}2-1}$. In short we can write $2^{\raise-1pt\hbox{$^{ \left\lfloor \frac{|P|+1}2 \right\rfloor}$}}$ for the dimension in both cases. Thus $HClS_I$ has dimension 
$2^{\raise-1pt\hbox{$^{n- \left\lfloor \frac{|P|+1}2 \right\rfloor}$}}$.

 A direct corollary to Theorem~\ref{MIEndo} is the following.

\begin{corollary}
  The induced module $M_I$ is the direct sum of
  $2^{\lfloor\frac{|P|+1}2\rfloor}$ isomorphic copies of $HClS_K$, where $K$ is
  the peak composition associated to $I$. \qed
\end{corollary}

We are now ready to show our main theorem and define the Frobenius characteristic between $HCl_n(0)$ modules and $\PP^*$, which we again denote by $\ch$.

\begin{theorem}
  The set $\{HClS_I :=Cl_n e_i \epsilon_I\}$, where $I$ runs over all peak
  compositions, is a complete set of pairwise non-isomorphic simple supermodules of
  $HCl_n(0)$. Moreover, there is a graded Hopf isomorphism defined by
  \begin{equation}
    \begin{array}{rccc}
      \ch\colon& \tilde{\mathcal G} & \longrightarrow &\PP^*   \\ \null \\
           &HClS_I & \longrightarrow & 
      \displaystyle  2^{\raise-1pt\hbox{$^{- \left\lfloor \frac{|P|+1}2 \right\rfloor}$}}  \Theta_{HP(I)}
            \end{array}
  \end{equation}
  where $HP(I)$ is the peak set of $I$, and $\tilde{\mathcal G}= \bigoplus_{n\geq 0} G_0(HCl_n(0))$.

\end{theorem}

Thus the $(1-q)$-transform at $q=-1$ can be interpreted as the induction map
from $ G_0(H_n(0)$ to $ G_0(HCl_n(0)$. This maps
${\mathcal G}$ to $\tilde{\mathcal G}$.

\medskip
\Proof Suppose that $S$ is a simple supermodule of $HCl_n(0)$. Decompose the
$H_n(0)$-socle of $S$ into simple modules
and choose a non-zero vector $v$ in
one of these simple factors. 
Then $v$ is a common eigenvector of all the $T_i$, 
so that there is a $I$ such that $\epsilon_I\mapsto v$ defines a
$H_n(0)$-morphism $\phi : S_I \to S$. 
Then, since $S$ is  supersimple, $v$
generates $S$ under the action of $HCl_n(0)$. Thus by induction, there is a
surjective morphism $M_I \to S$. Hence, each simple module of $HCl_n(0)$ must
be a quotient of some $M_I$ and consequently of some $HClS_I$.

Now, we know that given two $HClS_I$, either they are isomorphic (when they have the same peak sets) or else
there is no morphism between them. Thus $HClS_I$ has to be simple. The multiplication and comultiplication structures are induced from $\mathcal G$ and $QSym$ and the Frobenius characteristic between them.
\qed

By duality, we obtain $\ch^*\colon \PP\to \tilde{\mathcal G} ^*$ where $\tilde{\mathcal G} ^*= \bigoplus_{n\geq 0} K_0(HCl_n(0))$. We also remark that the dimension of the superradical of $HCl(0)$ is thus
  \begin{equation}
      2^nn!-\sum_P  2^{ 2n - (|P|+1) }\,.
  \end{equation}
where the sum is over all peak sets of $n$. It is still open to find nice generators for this radical.

\subsection{Decomposition matrices}
%%%%%%%%%%%%%%%%%%%%%%%%%%%%%%%%%%%

The generic Hecke-Clifford algebra $HCl_n(q)$ has its simple modules
$U_\lambda$ labelled by partitions into distinct parts,
and under restriction to $H_n(q)$, $U_\lambda$ has  as Frobenius
characteristic
\begin{equation}
\ch(U_\lambda)=2^{-\lfloor\ell(\lambda)/2\rfloor}Q_\lambda
\end{equation}
where $Q_\lambda$ is Schur's $Q$-function (see \cite{Ol,JN}).

By Stembridge's formula \cite{St},
we have
\begin{equation}
\ch(U_\lambda)=2^{-\lfloor\ell(\lambda)/2\rfloor}
\sum_{T\in{\mathcal ST}^\lambda}\Theta_{\Lambda(T)}
\end{equation}
where ${\mathcal ST}^\lambda$ is the set of standard shifted tableaux
of shape $\lambda$, and $\Lambda(T)$ the peak set of $T$.
The quasisymmetric characteristics of the simple $HCl_n(0)$
modules are proportional to the $\Theta$-functions, and the
coefficients $d_{\lambda I}$ in the
expression
\begin{equation}
\ch(U_\lambda)=\sum_I d_{\lambda I}\ch(HClS_I)
\end{equation}
form the decomposition matrices of the Hecke-Clifford algebras
at $q=0$.

For $\lambda$ a strict partition of $n$, and $I$ a peak composition of $n$
with peak set $P$, one has explicitly
\begin{equation}
d_{\lambda I}=
2^{\lfloor \frac12\ell(I)\rfloor -\lfloor\frac12\ell(\lambda)\rfloor}
\big| \{T\in{\mathcal ST}^\lambda\, |\, \Lambda(T)=P\} \big| \,.
\end{equation}
This is the analog for $HCl_n(0)$ of Carter's combinatorial formula
for the decomposition numbers of $H_n(0)$ \cite{Ca}.

\bigskip

Here are the  decomposition matrices for $n\le 9$. Note that
for $n=2,3$, $HCl_n(0)$ is semi-simple.
\begingroup
\tiny

\newlength{\matricelabelcellwidth}
\newlength{\matricelabelcellheight}
\newcolumntype{C}{@{\,}>{\hfil$}m{\matricelabelcellwidth}<{$}}
\newcolumntype{T}%
   {@{\,}>{\hfil\ \begin{rotate}{90}$}m{\matricelabelcellwidth}<{$\end{rotate}}}
\newcommand{\matricelabel}[6][21]{\vskip0.5cm%
  \settowidth{\matricelabelcellwidth}{#1}
  \settoheight{\matricelabelcellheight}{#1}
  \raisebox{1.33ex}%
  {$ %
    \begin{array}{c@{}c}
      & \begin{array}{#3}
        #4
      \end{array} \\
      \begin{array}{>{\vrule height \matricelabelcellheight width 0pt}r}
        #6
      \end{array}
      &
      \left[
        \begin{array}{#2}
          #5
        \end{array}
      \right] 
    \end{array}$}}

\matricelabel{CC}{TT}%
{3 & 21}%
{
  1 & 0\\
  0 & 1\\
}%
{3 \\ 21}%

\matricelabel{CCC}{TTT}%
{4 & 31 & 22}%
{
1 & 0 & 0\\
0 & 1 & 1\\
}%
{4 \\ 31}%

\matricelabel{CCCCC}{TTTTT}%
{3 & 41 & 32 & 23 & 221}%
{
                                1 & 0 & 0 & 0 & 0\\
                                0 & 1 & 1 & 1 & 0\\
                                0 & 0 & 1 & 0 & 1\\
}%
{5 \\ 41 \\ 32 }

\matricelabel{CCCCCCCC}{TTTTTTTT}%
{
  6 & 51 & 42 & 33 & 321 & 24 & 231 & 222
}%
{
1 & 0 & 0 & 0 & 0 & 0 & 0 & 0\\
0 & 1 & 1 & 1 & 0 & 1 & 0 & 0\\
0 & 0 & 1 & 1 & 1 & 0 & 1 & 1\\
0 & 0 & 0 & 0 & 1 & 0 & 0 & 1\\
}%
{ 6\\ 51\\ 42\\ 321 \\} 

\matricelabel{CCCCCCCCCCCCC}{TTTTTTTTTTTTT}
{ 7 & 61 & 52 & 43 & 421 & 34 & 331 & 322 & 25 & 241 & 232 & 223 & 2221 }%
{
  1 & 0 & 0 & 0 & 0 & 0 & 0 & 0 & 0 & 0 & 0 & 0 & 0\\
  0 & 1 & 1 & 1 & 0 & 1 & 0 & 0 & 1 & 0 & 0 & 0 & 0\\
  0 & 0 & 1 & 1 & 1 & 1 & 1 & 1 & 0 & 1 & 1 & 1 & 0\\
  0 & 0 & 0 & 1 & 0 & 0 & 1 & 1 & 0 & 0 & 1 & 0 & 2\\
  0 & 0 & 0 & 0 & 1 & 0 & 1 & 2 & 0 & 0 & 1 & 1 & 2\\
  }%
{ 7 \\ 61 \\ 52 \\ 43 \\ 421 \\  }%

\matricelabel{CCCCCCCCCCCCCCCCCCCCC}{TTTTTTTTTTTTTTTTTTTTT}
{ 8 & 71 & 62 & 53 & 521 & 44 & 431 & 422 & 35 & 341 & 332 & 323 & 3221 & 26 & 251 & 242 & 233 & 2321 & 224 & 2231 & 2222 }%
{
1 & 0 & 0 & 0 & 0 & 0 & 0 & 0 & 0 & 0 & 0 & 0 & 0 & 0 & 0 & 0 & 0 & 0 & 0 & 0 & 0\\
0 & 1 & 1 & 1 & 0 & 1 & 0 & 0 & 1 & 0 & 0 & 0 & 0 & 1 & 0 & 0 & 0 & 0 & 0 & 0 & 0\\
0 & 0 & 1 & 1 & 1 & 1 & 1 & 1 & 1 & 1 & 1 & 1 & 0 & 0 & 1 & 1 & 1 & 0 & 1 & 0 & 0\\
0 & 0 & 0 & 1 & 0 & 1 & 1 & 1 & 0 & 1 & 2 & 1 & 2 & 0 & 0 & 1 & 1 & 2 & 0 & 2 & 2\\
0 & 0 & 0 & 0 & 1 & 0 & 1 & 2 & 0 & 1 & 2 & 2 & 2 & 0 & 0 & 1 & 1 & 2 & 1 & 2 & 2\\
0 & 0 & 0 & 0 & 0 & 0 & 1 & 1 & 0 & 0 & 2 & 1 & 4 & 0 & 0 & 0 & 1 & 2 & 0 & 2 & 4\\
}%
{ 8 \\ 71 \\ 62 \\ 53 \\ 521 \\ 431 \\  }%

\matricelabel{CCCCCCCCCCCCCCCCCCCCCCCCCCCCCCCCCC}{TTTTTTTTTTTTTTTTTTTTTTTTTTTTTTTTTT}
{ 9 & 81 & 72 & 63 & 621 & 54 & 531 & 522 & 45 & 441 & 432 & 423 & 4221 & 36 & 351 & 342 & 333 & 3321 & 324 & 3231 & 3222 \
& 27 & 261 & 252 & 243 & 2421 & 234 & 2331 & 2322 & 225 & 2241 & 2232 & 2223 & 22221 }%
{
1 & 0 & 0 & 0 & 0 & 0 & 0 & 0 & 0 & 0 & 0 & 0 & 0 & 0 & 0 & 0 & 0 & 0 & 0 & 0 & 0 & 0 & 0 & 0 & 0 & 0 & 0 & 0 & 0 & 0 & 0 \
& 0 & 0 & 0\\
0 & 1 & 1 & 1 & 0 & 1 & 0 & 0 & 1 & 0 & 0 & 0 & 0 & 1 & 0 & 0 & 0 & 0 & 0 & 0 & 0 & 1 & 0 & 0 & 0 & 0 & 0 & 0 & 0 & 0 & 0 \
& 0 & 0 & 0\\
0 & 0 & 1 & 1 & 1 & 1 & 1 & 1 & 1 & 1 & 1 & 1 & 0 & 1 & 1 & 1 & 1 & 0 & 1 & 0 & 0 & 0 & 1 & 1 & 1 & 0 & 1 & 0 & 0 & 1 & 0 \
& 0 & 0 & 0\\
0 & 0 & 0 & 1 & 0 & 1 & 1 & 1 & 1 & 1 & 2 & 1 & 2 & 0 & 1 & 2 & 2 & 2 & 1 & 2 & 2 & 0 & 0 & 1 & 1 & 2 & 1 & 2 & 2 & 0 & 2 \
& 2 & 2 & 0\\
0 & 0 & 0 & 0 & 1 & 0 & 1 & 2 & 0 & 1 & 2 & 2 & 2 & 0 & 1 & 2 & 2 & 2 & 2 & 2 & 2 & 0 & 0 & 1 & 1 & 2 & 1 & 2 & 2 & 1 & 2 \
& 2 & 2 & 0\\
0 & 0 & 0 & 0 & 0 & 1 & 0 & 0 & 0 & 1 & 1 & 1 & 0 & 0 & 0 & 1 & 1 & 2 & 0 & 2 & 2 & 0 & 0 & 0 & 1 & 0 & 0 & 2 & 2 & 0 & 0 \
& 2 & 0 & 2\\
0 & 0 & 0 & 0 & 0 & 0 & 1 & 1 & 0 & 1 & 3 & 2 & 4 & 0 & 0 & 2 & 3 & 6 & 1 & 6 & 8 & 0 & 0 & 0 & 1 & 2 & 1 & 4 & 6 & 0 & 2 \
& 6 & 4 & 4\\
0 & 0 & 0 & 0 & 0 & 0 & 0 & 0 & 0 & 0 & 1 & 0 & 2 & 0 & 0 & 0 & 1 & 2 & 0 & 2 & 4 & 0 & 0 & 0 & 0 & 0 & 0 & 2 & 2 & 0 & 0 \
& 2 & 2 & 2\\
}%
{ 9 \\ 81 \\ 72 \\ 63 \\ 621 \\ 54 \\ 531 \\ 432 \\  }%

\endgroup

\end{document}